\documentclass[final,leqno,onefignum,onetabnum]{siamltex}

\title{Spectral deferred corrections with fast-wave slow-wave splitting}

\author{Daniel Ruprecht\footnotemark[1] \and Robert Speck\footnotemark[2]}

\usepackage{amsmath}
\usepackage{amssymb}
\usepackage{color}
\usepackage{graphicx}
\usepackage[free-standing-units,group-separator={,}]{siunitx}
\usepackage[colorlinks=true]{hyperref}

\newcommand{\Tt}[1]{\mathbf{#1}}
\newcommand{\ffast}[0]{ f_{\text{f}} }
\newcommand{\fslow}[0]{ f_{\text{s}} }
\newcommand{\lamfast}[0]{ \lambda_{\text{f}} }
\newcommand{\lamslow}[0]{ \lambda_{\text{s}} }
\newcommand{\Qfast}[0]{\Tt{Q}_{\Delta}^{\rm fast}}
\newcommand{\Qslow}[0]{\Tt{Q}_{\Delta}^{\rm slow}}

\newcommand{\fwsw}{\textsc{fwsw}-SDC}
\newtheorem{remark}{Remark}

\begin{document}

\maketitle

\begin{abstract}
The paper investigates a variant of semi-implicit spectral deferred corrections (SISDC) in which the stiff, fast dynamics correspond to fast propagating waves (``fast-wave slow-wave problem'').
We show that for a scalar test problem with two imaginary eigenvalues $i \lamfast$, $i \lamslow$, having $\Delta t \left( \left| \lamfast \right| + \left| \lamslow \right| \right) < 1$ is sufficient for the fast-wave slow-wave SDC (\fwsw) iteration to converge and that in the limit of infinitely fast waves the convergence rate of the non-split version is retained.
Stability function and discrete dispersion relation are derived and show that the method is stable for essentially arbitrary fast-wave CFL numbers as long as the slow dynamics are  resolved.
The method causes little numerical diffusion and its semi-discrete phase speed is accurate also for large wave number modes.
Performance is studied for an acoustic-advection problem and for the linearised Boussinesq equations, describing compressible, stratified flow.
\fwsw~is compared to a diagonally implicit Runge-Kutta (DIRK) and IMEX Runge-Kutta (IMEX) method and found to be competitive in terms of both accuracy and cost.
\end{abstract}


\renewcommand{\thefootnote}{\fnsymbol{footnote}}

\footnotetext[1]{School of Mechanical Engineering, University of Leeds, Leeds LS2 9JT, UK}
\footnotetext[2]{J\"ulich Supercomputing Centre, Forschungszentrum J\"ulich GmbH, Germany}

\slugger{mms}{xxxx}{xx}{x}{x--x}

\renewcommand{\thefootnote}{\arabic{footnote}}

\pagestyle{myheadings}
\thispagestyle{plain}
\markboth{D. Ruprecht and R. Speck}{Spectral deferred corrections with fast-wave slow-wave splitting} 


\section{Introduction}
For simulations of compressible flow, in particular in numerical weather prediction and climate simulations, the presence of acoustic waves can pose significant numerical challenges to the time integration method.
Explicit methods are restricted to inefficiently small steps while fully implicit methods are expensive and can artificially slow down high wave number modes.
The fully compressible equations can be replaced by filtered models that do not support sound waves~\cite{Durran1989,OguraPhillips1962}, but these require solution of a Poisson problem in each step and have difficulties capturing large-scale wave dynamics~\cite{DaviesEtAl2003}.

Therefore, a widely used class of methods are split-explicit integrators: they separate the equation into fast and slow processes which are then integrated with different time step sizes and different (explicit) methods.
A popular method of this type is a third-order Runge-Kutta scheme combined with a forward-backward Euler integrator for the acoustic terms~\cite{WickerSkamarock2002}.
While computationally efficient, split-explicit methods typically require some form of damping for stabilisation~\cite{Baldauf2010}, which reduces their effective order of accuracy.
However, a second-order split-explicit two-step peer method has recently been derived that allows for stable integration of the compressible Euler equations without damping~\cite{KnothEtAl2014}.

Another form of splitting are semi-implicit methods.
They also split the equations into fast and slow parts but then use an implicit method for the fast and an explicit method for the slow part.
In many applications, the fast, stiff terms stem from diffusion and/or rapid chemical reactions and methods with IMEX splitting for equations of reaction-diffusion type have been widely studied~\cite{Koto2008, NieEtAl2006, ShampineEtAl2006}.
For stratified, compressible flows, however, the fast dynamics are not diffusive but stem from acoustic and fast gravity waves while the slow dynamics correspond to slower waves and advection.
IMEX splitting methods for such ``fast-wave slow-wave" problems~\cite{DurranEtAl2012} have not been as widely studied. 
Some literature does exist~\cite{GiraldoEtAl2010,Weller2013}, however, and early works go back to the 1970's~\cite{Kwizak1971,TappWhite1976}.
The performance of IMEX Runge-Kutta methods has only recently been studied for fast-wave slow-wave problems~\cite{WhitakerKar2013}, inspired by a previous study for multi-step methods of IMEX-type~\cite{DurranEtAl2012}.
A general framework for both multi-step and Runge-Kutta IMEX methods for the ``Nonhydrostatic  Unified Model of the Atmosphere" has been recently developed and tests found that higher-order time stepping methods are more efficient~\cite{GiraldoEtAl2013}. 
Splitting methods for use in climate simulations are also an active topic of research~\cite{CollinsEtAl2015}.

Derivation of high-order IMEX methods can be difficult and leads to a quickly growing number of order conditions~\cite{PareschiRusso2005}.
Third-order four-stage IMEX methods have been derived~\cite{AscherEtAl1997,PareschiRusso2005} as well as a fourth-order method with six stages and a fifth-order method with eight stages~\cite{KennedyCarpenter2003}.
In contrast, semi-implicit spectral deferred corrections (SISDC)~\cite{Minion2003} allow for the simple and generic construction of split methods of arbitrary order.
SISDC have been studied and found to be competitive for advection-reaction-diffusion problems~\cite{BourliouxEtAl2003,LaytonMinion2004}.
Also, it has been shown that, for smooth solutions and Lipshitz continuous right-hand sides, SISDC can attain the full accuracy of the underlying collocation formula~\cite{HagstromZhou2006}.
Defect correction methods with splitting based on equidistant instead of spectral nodes have also recently been investigated~\cite{Christlieb2015}.
However, the performance of SISDC for fast-wave slow-wave problems has only been analysed rudimentarily so far~\cite{Weingartz2014}.

This paper investigates the performances of SISDC with ``fast-wave slow-wave'' splitting (\fwsw).
Convergence of \fwsw~is shown for the case where both wave types are well resolved and in the limit of infinitely fast acoustic waves.
We derive the stability function of \fwsw~and show that the method possesses favourable stability characteristics: for a reasonable range of slow wave speeds, the method remains stable for arbitrarily large fast wave speeds.
The semi-discrete dispersion relation is derived and shows that \fwsw~damps high wave number modes (which are typically spatially under-resolved) while correctly propagating other modes.
Finally, the iterative nature of SDC produces increasingly accurate starting values for whatever iterative solver is used for the implicit part.
We demonstrate that \fwsw~can be more efficient than a diagonally implicit Runge-Kutta method (DIRK) of the same order~\cite{Alexander1977,KennedyCarpenter2016} and that it can compete with Runge-Kutta IMEX methods: even though SDC requires significantly more linear systems to be solved, the total number of required GMRES iterations is only slightly larger (or even comparable) because the increasingly accurate starting values lead to rapid convergence. 

\section{Spectral deferred corrections}
Consider an initial value problem of the following form
\begin{equation}
	\label{eq:ivp}
	u'(t) = f(u(t)), \quad u(t_0) = u_0.
\end{equation}
For the sake of simplicity, we consider integration of~\eqref{eq:ivp} over one time step $[T_n, T_{n+1}]$ with length $\Delta t := T_{n+1} - T_n$. 
We also focus on the autonomous case, but the extension to the non-autonomous vector case is straightforward.

\subsection{Collocation}
For smooth solutions, the initial value problem in differential form~\eqref{eq:ivp} is equivalent to the integral equation
\begin{equation}
	\label{eq:integral}
	u(t) = u(T_0) + \int_{T_0}^{t} f(u(s))~\mathrm{d}s, \quad T_n \leq t \leq T_{n+1}.
\end{equation}
We introduce $M$ quadrature nodes\footnote{Throughout the paper we consider Radau nodes, see also the comments in Section~\ref{subsec:quad_node_choice}.} $T_n \leq \tau_1 < \ldots < \tau_M \leq T_{n+1}$ and denote as $\Delta \tau_m := \tau_{m} - \tau_{m-1}$ for $m=2, \ldots, M$ the distance between two nodes.
For $m=1$, we define $\Delta \tau_1 := \tau_1 - T_n$.
Note that for nodes where $\tau_1$ coincides with $T_n$ (e.g. Gauss-Lobatto nodes) we have $\Delta \tau_1 = 0$.
We approximate the integral in~\eqref{eq:integral} by the corresponding quadrature rule to get the collocation equations
\begin{equation}
	\label{eq:picard}
	u_m = u_0 + \sum_{j=1}^{M} q_{m,j} f(u_j), \quad m=1,\ldots,M.
\end{equation}
Here, $u_0 \approx u(T_n)$ is the initial value brought forward from the previous time step, $u_m \approx u(\tau_m)$ is the approximate solution at quadrature point $\tau_m$ while the $q_{m,j}$ are weights defined as
\begin{equation}
	\label{eq:weights}
	q_{m,j} := \int_{T_n}^{\tau_m} l_j(s)~\mathrm{d}s, \quad m,j=1,\ldots,M,
\end{equation}
with $l_j$ being the Lagrange polynomials to the points $\tau_m$.
Once the stages $u_j$ are known, the final update step
\begin{equation}
	\label{eq:picard_update}
	u_{n+1} = u_0 + \sum_{j=1}^{M} q_j f(u_j)
\end{equation}
provides $u_{n+1} \approx u(T_{n+1})$ with
\begin{equation}
	q_j := \int_{T_n}^{T_{n+1}} l_j(s)~\mathrm{d}s, \quad j=1,\ldots,M.
\end{equation}
Solving~\eqref{eq:picard} for the stages directly and using the $u_m$ in the update~\eqref{eq:picard_update} corresponds to a collocation method.
Collocation methods are a subclass of implicit Runge-Kutta methods with the $u_m$ being the stages and the $q_{m,j}$ the entries in the Butcher tableau~\cite[Theorem 7.7]{HairerEtAl1993_nonstiff}.
They require solving one large system composed of the $M$ coupled nonlinear equations~\eqref{eq:picard}.
\begin{remark}
By using weights $\tilde{q}_j := \int_{T_n}^{\theta} l_j(s)~\mathrm{d}s$ in~\eqref{eq:picard_update}, an $M^{\text{th}}$-order accurate approximate solution can be constructed at any value $T_n \leq \theta \leq T_{n+1}$, thereby naturally providing a dense output~\cite[Sect. II.6]{HairerEtAl1993_nonstiff} formula.
\end{remark}
\begin{remark}\label{remark:update_lobatto}
For nodes where $\tau_1 = T_{n}$, we get $q_{1,j} = 0$ for $j=1,\ldots,M$ from~\eqref{eq:weights} so that~\eqref{eq:picard} for $m=1$ reduces to $u_1 = u_0$.
Analogously, if $\tau_M = T_{n+1}$, we have $q_{M,j} = q_{j}$ for $j=1,\ldots,M$ and~\eqref{eq:picard} for $m=M$ is identical to~\eqref{eq:picard_update} so that $u_{n+1} = u_M$.
\end{remark}

\subsection{Spectral deferred corrections}
Instead of directly solving for the intermediate solutions $u_m$, spectral deferred corrections (SDC)~\cite{DuttEtAl2000} proceed with the following iteration that avoids solving the fully coupled system~\eqref{eq:picard} and solves a series of smaller problems instead.
With implicit Euler as base method, the SDC iteration reads
\begin{equation}
	\label{eq:implicit_sdc_sweep}
	u^{k+1}_{m} = u^{k+1}_{m-1} + \Delta \tau_m \left( f(u^{k+1}_{m}) - f(u^k_m) \right) + \sum_{j=1}^{M} s_{m,j} f(u^k_m), \quad m=1,\ldots,M
\end{equation}
with $u^k_0 = u_0$, $s_{m,j} := q_{m,j} - q_{m-1,j}$ for $m=2,\ldots,M$, $s_{1,j} := q_{1,j}$ and $k$ being the iteration index.
\begin{remark}
If iteration~\eqref{eq:implicit_sdc_sweep} converges and $u^{k+1}_m - u^{k}_m \to 0$, it reduces to
\begin{equation}
	u_m = u_{m-1} + \sum_{j=1}^{M} s_{m,j} f(u_m)
\end{equation}
from which it readily follows that
\begin{equation}
	u_m = u_0 + \sum_{l=1}^{m} \sum_{j=1}^{M} s_{l,j} f(u_j) = u_0 + \sum_{j=1}^{M} q_{m,j} f(u_j).
\end{equation}
Therefore, if SDC converges it reproduces the collocation solution~\eqref{eq:picard} at each $\tau_m$.
\end{remark}

For the scalar case, we later derive an upper bound for the convergence rate for small enough $\Delta t$.
However, the attractiveness of SDC stems from the fact that full convergence is not required to produce a useful approximation of $u(T_{n+1})$.
It has been shown that using $k$ iterations with either implicit or explicit Euler as base method results in a $k^{\text{th}}$-order method if the underlying quadrature rule is sufficiently accurate~\cite{ShuEtAl2007}.
Higher order methods can be used as SDC base method but do not necessarily improve the order by more than one per iteration~\cite{ChristliebEtAl2010_MoC}.
SDC can also be written as a preconditioned iteration for the solution of~\eqref{eq:picard}~\cite{HuangEtAl2006}.
For approximate stages $\tilde{u}_m \approx u_m$, the components of the residual are defined as
\begin{equation}
	\label{eq:residual}
	r_m =  u_0 + \sum_{j=1}^{M} q_{m,j} f(\tilde{u}_j) - \tilde{u}_m
\end{equation}
and can be used to monitor convergence.
This interpretation has been used to derive a number of modifications of SDC~\cite{SpeckEtAl2015_BIT,SpeckEtAl2015_DDM}.
SDC can also be used as framework for the derivation of high-order multi-rate methods~\cite{BourliouxEtAl2003,EmmettEtAl2014}.

\subsection{Semi-implicit SDC}
Consider now a case where the right-hand side of the initial value problem~\eqref{eq:ivp} can be split into a fast and a slow term as
\begin{equation}
		u'(t) = f(u(t)) = \ffast(u(t)) + \fslow(u(t)), \quad u(t_0) = u_0.
\end{equation}
Typically, $\ffast$ and $\fslow$ come from the spatial discretisation of different terms of a partial differential equation.
IMEX Euler can be used as base method, treating the slow part explicitly and the fast part implicitly.
In this case, the SDC iteration~\eqref{eq:implicit_sdc_sweep} becomes
\begin{equation}
	u^{k+1}_m = u^{k+1}_{m-1} + \Delta \tau_m \left( \ffast(u^{k+1}_m) - \ffast(u^k_m) + \fslow(u^{k+1}_{m-1}) - \fslow(u^k_{m-1})\right) + \sum_{j=1}^{M} s_{m,j} f(u^k_m)
\end{equation}
for $m=1,\ldots,M$.
Previous works have analysed the case where $\ffast$ is a term describing diffusion or a fast chemical reaction~\cite{BourliouxEtAl2003,LaytonMinion2004,Minion2003}.
Here, we analyse purely hyperbolic problems in which both $\ffast$ and $\fslow$ stem from the discretisation of terms describing wave propagation but at different speeds (``fast-wave slow-wave SDC'' or \fwsw~for short).
An important example are atmospheric flows, where ``physically insignificant fast waves''~\cite[Chap. 8]{Durran2010} like acoustic and fast gravity waves impose severe limitations on time steps for explicit methods compared to e.g.~slow moving Rossby waves or advection.

\section{Theory}
There are two different vantage points from which \fwsw~can be analysed: as a split method with a fixed order set by a fixed number of iterations $K$ for a sufficiently large number of nodes $M$, or as an iterative solver for the collocation problem where iterations are performed until the norm of the residual~\eqref{eq:residual} reaches a prescribed tolerance.
We investigate \fwsw~from both viewpoints for the scalar test problem
\begin{equation}
	\label{eq:scalar}
	u_t(t) = i \lambda u(t) = i \lamfast u(t) + i \lamslow u(t), \quad u(0) = 1, \quad \lamfast, \lamslow \in \mathbb{R}
\end{equation}
with $\lamfast \gg \lamslow$. 
Convergence towards the collocation solution is assessed by analysing norm and spectral radius of the error propagation matrix.
Then, for a fixed number of iterations, stability is analysed by deriving the stability function for \fwsw.
Finally, also for fixed $K$, the semi-discrete dispersion relation for \fwsw~applied to an acoustic-advection problem is derived and wave propagation characteristics are analysed.

Model problem~\eqref{eq:scalar} and the term ``fast-wave slow-wave" are borrowed from recent work analysing multi-step methods with IMEX splitting~\cite{DurranEtAl2012}.
Equation~\eqref{eq:scalar} is frequently used to investigate stability of integration schemes for meteorological applications~\cite{DurranEtAl2012,Weller2013}.
Note that, in contrast to the standard Dahlquist test equation,~\eqref{eq:scalar} features no real eigenvalue but two imaginary eigenvalues of different magnitude.
When applied to~\eqref{eq:scalar} the SDC sweep~\eqref{eq:implicit_sdc_sweep} becomes
\begin{equation}
	\label{eq:fwsw-sdc-sweep}
	u^{k+1}_{m} = u^{k+1}_{m-1} + \Delta \tau_m \left[ i \lamfast \left( u^{k+1}_m - u^k_m \right) +  i \lamslow \left( u^{k+1}_{m-1} - u^k_{m-1} \right) \right]
	+ \sum_{j=1}^{M} s_{m,j} i \lambda u^k_j.
\end{equation}
By recursively using~\eqref{eq:fwsw-sdc-sweep} it is straightforward to show that this ``node-to-node" formulation -- updating from $u_{m-1}$ to $u_m$ -- is equivalent to the ``zero-to-node" formulation
\begin{equation}
	\label{eq:fwsw-sdc-zero-to-node-sweep}
	u^{k+1}_m = u_0 + \sum_{j=1}^{m} \Delta \tau_j \left[ i \lamfast \left( u^{k+1}_{j} - u^k_j \right) +  i \lamslow \left( u^{k+1}_{j-1} - u^{k}_{j-1} \right) \right] + \sum_{j=1}^{M} q_{m,j} i \lambda u^k_j
\end{equation}
updating from $u_0$ to $u_m$ with the $q_{m,j}$ defined according to~\eqref{eq:weights}.
Collecting all intermediate solutions (i.e.~the stages) in a vector
\begin{equation}
	\Tt{U}^k := \left( u^k_1, \ldots, u^k_{M} \right)
\end{equation}
allows to compactly write~\eqref{eq:fwsw-sdc-zero-to-node-sweep} as
\begin{equation}
	\Tt{U}^{k+1} = \Tt{U}_0 + \Delta t \left[ \Tt{Q}^{\text{fast}}_{\Delta} i \lamfast \left( \Tt{U}^{k+1} - \Tt{U}^k \right) + \Tt{Q}^{\text{slow}}_{\Delta} i \lamslow \left( \Tt{U}^{k+1} - \Tt{U}^{k} \right) \right] + \Delta t \Tt{Q} i \lambda \Tt{U}^{k}
\end{equation}
with matrices
\begin{equation}
	\Qfast :=  \frac{1}{\Delta t}\begin{pmatrix} \Delta \tau_1 & \\ \Delta \tau_1 & \Delta \tau_2 & \\ \vdots & \vdots & \\ \Delta \tau_1 & \Delta \tau_2 & \ldots & \Delta \tau_M \end{pmatrix}
\end{equation}
and 
\begin{equation}
    \Qslow := \frac{1}{\Delta t}\begin{pmatrix} 0 & \\ \Delta \tau_1& 0 \\ \Delta \tau_1 & \Delta \tau_2 & 0 \\ \vdots & \vdots & \\ \Delta \tau_1 & \Delta \tau_2 & \ldots & \Delta \tau_{M-1} & 0 \end{pmatrix}
\end{equation}
and $\Tt{Q} = \left( q_{m,j}/\Delta t \right)_{m,j=1,\ldots,M}$ and $\Tt{U}_0 = \left( u_0, \ldots, u_0 \right)$.
Rearranging terms gives
\begin{equation}
	\label{eq:scalar-zero-to-node-system}
	\left( \Tt{I} - \Delta t \left( i \lamfast \Qfast + i \lamslow \Qslow \right) \right) \Tt{U}^{k+1} = \Tt{U}_0 + \Delta t \left( i \lambda \Tt{Q}  - \left( i \lamfast \Qfast +i  \lamslow \Qslow \right)\right) \Tt{U}^k.
\end{equation}
This is the \fwsw~iteration written as a preconditioned Richardson iteration to solve the collocation equation
\begin{equation}
	\label{eq:collocation}
	\Tt{U} = \Tt{U}_0 + \Delta t i \lambda \Tt{Q} \Tt{U}.
\end{equation}
An interesting variant of SDC (colloquially known as ``St.~Martin's trick") uses a LU decomposition instead of the above $\Tt{Q}_{\Delta}$'s (in particular instead of $\Qfast$) as a preconditioner~\cite{Weiser2014}. 
Investigating how this strategy affects \fwsw~is left for future work.

\subsection{Iteration error and local truncation error}
Because the solution $\Tt{U}$ of the collocation equation~\eqref{eq:collocation} is a fixed point of~\eqref{eq:scalar-zero-to-node-system}, the error $\Tt{e}^k := \Tt{U}^k - \Tt{U}$ between the exact collocation solution and its approximation $\Tt{U}^k$ provided by SDC after $k$ sweeps propagates according to
\begin{align}
	\Tt{e}^{k+1} &=  \left( \Tt{I} - \Delta t \left( i \lamfast \Qfast + i \lamslow \Qslow \right) \right)^{-1} \Delta t \left( i \lambda \Tt{Q} - \left( i \lamfast \Qfast + i \lamslow \Qslow \right) \right)  \Tt{e}^{k} \nonumber \\ &=: \Tt{E} \Tt{e}^{k}
	\label{eq:error_matrix}
\end{align}
with $\Tt{e}^0 = \Tt{U} - \Tt{U}_0$.
Below, we will derive a bound for the norm of the error propagation matrix $\Tt{E}$.
Using this bound we can show that \fwsw~converges and increases the order by one per iteration, up to the order of the collocation formula, if $\Delta t \left( \left| \lamfast \right| + \left| \lamslow \right| \right) < 1 $ (``non-stiff case'').
We also compute numerically the spectral radius of $\Tt{E}$ in the limit $\lamfast \to \infty$ (``stiff limit'') and show that it remains smaller than unity.
Therefore, \fwsw~also converges for $k \to \infty$ in the limit of infinitely fast acoustic waves as long as $\Delta t\left|  \lamslow \right|$ is small enough.
Moreover, as shown in Section~\ref{subsec:stability}, \fwsw~remains stable for arbitrary large values of $\lamfast$ even for a fixed small number of iterations if $\Delta t \left|  \lamslow \right|$ is small enough.
\subsubsection{Non-stiff case}
For the case where $\Delta t \left( \left| \lamfast \right| + \left| \lamslow \right| \right) < 1$ we give a simple proof that the \fwsw~iteration, using a combination of forward and backward Euler, converges and that each iteration increases the order by one.
A qualitative proof along similar lines (using a Neumann series expansion of the iteration matrix) for \textit{either} backward \textit{or} forward Euler method as base integrator has been given before~\cite{HuangEtAl2006}.
A proof for the generic case with splitting is also available~\cite{HagstromZhou2006}, but is more involved and does not directly provide an estimate for the iteration error.
\begin{lemma}\label{lemma:q_delta_bound}
For any set of quadrature nodes $(\tau_m)_{m=1,\ldots,M}$ in $[T_n,T_{n+1}]$ we have
\begin{equation}
	\left\| \Qfast \right\|_{\infty} \leq 1 \quad \text{and} \quad \left\| \Qslow \right\|_{\infty} \leq 1.
\end{equation}
\begin{proof} Since $\sum_{j=1}^{M} \Delta \tau_j \leq \Delta t$ it holds that
\begin{equation}
	\left\| \Qfast \right\|_{\infty} = \max_{i=1,\ldots,M}  \Delta t^{-1} \sum_{j=1}^{i} \Delta \tau_j \leq \Delta t^{-1} \Delta t = 1
\end{equation}
and analogously for $\Qslow$.
\end{proof}
\end{lemma}
\begin{lemma}\label{lemma:precond_bound}
If the time step $\Delta t$ is small enough so that $\Delta t \left( \left| \lamfast \right| + \left| \lamslow \right| \right) < 1$, it holds that 
\begin{equation}
	\left\| \left( \Tt{I} - \Delta t \left( i \lamfast \Qfast + i \lamslow \Qslow \right) \right)^{-1}\right\|_{\infty} \leq 1 + \Delta t \left( \left| \lamfast \right| + \left| \lamslow \right| \right)+ \mathcal{O}(\Delta t^2)
\end{equation}
\begin{proof}
By Lemma~\ref{lemma:q_delta_bound} we have
\begin{equation}
	\Delta t \left\| i \lamfast \Qfast + i \lamslow \Qslow  \right\|_{\infty} \leq \Delta t \left( \left| \lamfast \right| + \left| \lamslow \right| \right).
\end{equation} 
Therefore, if $\Delta t  \left( \left| \lamfast \right| + \left| \lamslow \right| \right) < 1$, the inverse matrix can be expanded in a Neumann series
\begin{equation}
	\left( \Tt{I} - \Delta t \left( i \lamfast \Qfast + i \lamslow \Qslow \right) \right)^{-1} = \Tt{I} + \Delta t \left( i \lamfast \Qfast + i \lamslow \Qslow \right) + \ldots
\end{equation}
Taking the norm and using Lemma~\ref{lemma:q_delta_bound} again shows the estimate.
\end{proof}
\end{lemma}
\begin{lemma}\label{lemma:q_bound}
For any set of quadrature nodes $(\tau_m)_{m=1,\ldots,M}$ in $[T_n, T_{n+1}]$ it holds that
\begin{equation}
	\left\| \Tt{Q} \right\|_{\infty} \leq \Lambda_M
\end{equation}
where
\begin{equation}
	\Lambda_M := \max_{-1 \leq x \leq 1} \sum_{j=1}^{M} \left| \tilde{l}_j(x) \right|
\end{equation}
is the Lebesgue constant and $\tilde{l}_j$ are the Lagrange polynomials on the interval $[-1,1]$.
\begin{proof}
We can transform the Lagrange polynomials on $[T_n, T_{n+1}]$ to $[-1,1]$ via the transformation 
\begin{equation}
	t \mapsto x = 2\frac{t - T_n}{\Delta t} - 1  \quad \text{with inverse} \quad x \mapsto t  = \left( \frac{x+1}{2}  \right) \Delta t + T_n.
\end{equation}
Therefore, by substitution,
\begin{equation}
	\left| q_{m,j} \right|  = \frac{1}{\Delta t} \left| \int_{T_n}^{\tau_m} l_j(s)~ds \right| \leq \frac{1}{\Delta t} \int_{T_n}^{T_{n+1}} \left| l_j(s) \right|~ds = \frac{1}{2} \int_{-1}^{1} \left| \tilde{l}_j(x) \right|~dx.
\end{equation}
Now we can compute
\begin{equation}
	\left\| \Tt{Q} \right\|_{\infty} = \max_{m=1,\ldots,M} \sum_{j=1}^{M} \left| q_{m,j} \right| \leq \frac{1}{2} \int_{-1}^{1} \sum_{j=1}^{M} \left| \tilde{l}_{j}(x) \right|~dx \leq \Lambda_M.
\end{equation}
\end{proof}
\end{lemma}

The following theorem is now readily proven:
\begin{theorem}\label{theorem:error_matrix_bound}
For $\Delta t \left( \left| \lamfast \right| + \left| \lamslow \right| \right) < 1 $, the norm of the error propagation matrix $\Tt{E}$ is bounded by
\begin{equation}
	\left\| \Tt{E} \right\|_{\infty} \leq \Delta t \left( \Lambda_M + \left| \lamfast \right| + \left| \lamslow \right| \right) + \mathcal{O}(\Delta t^2).
\end{equation}
\begin{proof}
Follows directly from Lemmas~\ref{lemma:precond_bound} and~\ref{lemma:q_bound}
\end{proof}
\end{theorem}

Since $\Lambda_M$, $\lamfast$ and $\lamslow$ are all independent of $\Delta t$, this estimate guarantees that \fwsw~eventually converges if $\Delta t$ becomes small enough and $\left\| \Tt{E} \right\|_{\infty} < 1$.
However, this condition is sufficient but not necessary and typically SDC already converges for time steps much larger than what could be expected from Theorem~\ref{theorem:error_matrix_bound}.
In particular, as shown below, \fwsw~converges and remains stable for arbitrarily large values of $\lamfast$.
Also, the provided bound is not sharp.
One reason seems to be that Lemma~\ref{lemma:q_bound} gives a very pessimistic estimate of the norm  of $\Tt{Q}$, at least for spectral nodes.
Numerical experiments not documented here suggest that actually $\left\| \Tt{Q} \right\|_{\infty} \leq 1$ might hold for Lobatto, Radau and Legendre nodes, but we do not have a rigorous proof for this hypothesis.
In addition, it may be more favourable to estimate the norm of the difference between $\lambda\Tt{Q}$ and $\lamfast \Qfast +  \lamslow \Qslow$ in~\eqref{eq:error_matrix}, but a promising approach to do this has not yet been found.

For the case where both fast and slow waves are well resolved, we can now show that the local truncation error of \fwsw~with $k$ iterations is of order $k+1$, up to the order of the underlying quadrature rule.
Assume that $u_0 = u(T_n)$ is the exact solution at the beginning of the time step.
Denote as $u_{n+1}$ the solution at the end of the time step generated by~\eqref{eq:picard_update} using the exact stages of the collocation solution $u_m$.
Further denote as $u^k_{n+1}$ the solution also computed from~\eqref{eq:picard_update} but using the approximate stages $u^k_m$ computed with $k$ sweeps of \fwsw.
Then,
\begin{equation}
	u_{n+1} - u^k_{n+1} = i \lambda \sum_{j=1}^{M} q_{j} \left( u_j - u^k_j \right).
\end{equation}
According to Theorem~\ref{theorem:error_matrix_bound}, the difference between the exact stages $u_m$ and the approximate stages $u_m^k$ satisfies 
\begin{equation}
\left\| \Tt{e}^{k} \right\|_{\infty} = \left\| \Tt{U} - \Tt{U}^{k} \right\|_{\infty} = \max_{m=1,\ldots,M} \left| u_m - u^k_m \right| = \mathcal{O}(\Delta t^k).
\end{equation}
Also, by a similar argument as in the proof of Lemma~\ref{lemma:q_bound}, we have $\left| q_j \right| = \mathcal{O}(\Delta t)$ for all $j=1,\ldots,M$.
Together, this gives
\begin{equation}
	\left| u_{n+1} - u_{n+1}^k \right| \leq \left| \lambda \right| \sum_{j=1}^{M} \left| q_j \right| \left| u_j - u^k_j \right| = \mathcal{O}(\Delta t^{k+1}).
\end{equation}
For the collocation solution, that is the $u_m$ which satisfy~\eqref{eq:picard} exactly, the truncation error at the end of the step is
\begin{equation}
	\left| u(T_{n+1}) - u_{n+1} \right| = \left| \lambda \right| \left| \int_{T_n}^{T_{n+1}} u(s)~ds -\sum_{j=1}^{M} q_{j} u_j \right|= \mathcal{O}(\Delta t^{p+1})
\end{equation}
where $p$ is the order of the quadrature rule.
For Lobatto nodes we would have $p = 2M-2$, for Radau nodes $p=2M-1$ and for Legendre nodes $p=2M$.
The local truncation error of \fwsw~thus is, using triangle inequality,
\begin{equation}
	\label{eq:lte}
	\left| u(T_{n+1}) - u^k_{n+1} \right| = \mathcal{O}(\Delta t^{k+1}) + \mathcal{O}(\Delta t^{p+1}) = \mathcal{O}(\Delta t^{ \min\left\{ k+1, p+1 \right\} } ).
\end{equation}
The same result has been previously derived for SDC with \textit{either} implicit \textit{or} explicit Euler as base method using a different approach based on induction~\cite{ShuEtAl2007}.
While the proof can be adopted for \fwsw, the  approach presented here provides an explicit estimate for the \emph{iteration error}, that is the difference between the SDC and the collocation solution.
This is beneficial when SDC is not used to generate a method with a fixed order but iterations are instead performed until some residual tolerance is reached.
Also, the interpretation and analysis of SDC as a linear iteration can provide a starting point for the mathematical analysis of SDC's multi-level variants MLSDC and PFASST. Such an analysis will be pursued in future work.
\begin{remark}\label{remark:update}
When $T_{n+1}$ is a quadrature node (e.g. for Gauss-Lobatto nodes), one can simply set $u_{n+1} = u_M$ instead of performing update~\eqref{eq:picard_update}.
For the exact collocation solution this makes no difference (see Remark~\ref{remark:update_lobatto}) but if the stages are only approximately computed then the two updates give different results.
Experiments not documented here suggest that setting $u_{n+1} = u_M$ gives a slightly less accurate approximation but can significantly improve stability for Gauss-Lobatto nodes and might therefore be a useful strategy.
\end{remark}

\subsubsection{Stiff limit}
\begin{figure}[th!]
	\centering	
	\includegraphics[scale=1]{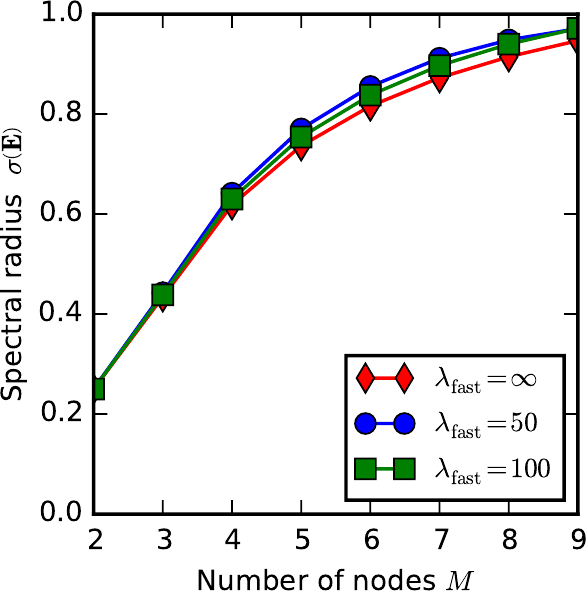}\hfill
	\includegraphics[scale=1]{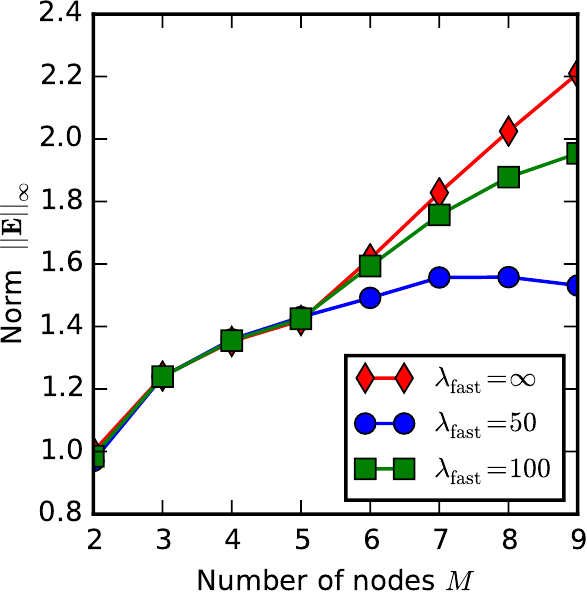}
	\caption{Spectral radius (left) and norm (right) of the error propagation matrix $\Tt{E}$ in the limit $\lamfast \to \infty$ (red) and for large but finite values of $\lamfast = 50$ (blue) and $\lamfast = 100$ (green). All cases use Gauss-Radau nodes, $\Delta t = 1.0$ and $\lamslow = 1.0$.}
	\label{fig:stiff_limit}
\end{figure}

One key advantage of \fwsw~is that the splitting does not impair convergence: even in the limit of infinitely fast fast waves, \fwsw~converges as good (or bad) as the non-split version based on backward Euler.
For fixed $\Delta t$ and $\lamslow$, in the limit $\lamfast \to \infty$ the error propagation matrix~\eqref{eq:error_matrix} becomes
\begin{equation}
	\label{eq:error_matrix_limit}
	\Tt{E} = \Tt{I} - (\Qfast)^{-1} \Tt{Q}.
\end{equation}
This is identical to the stiff limit of non-split SDC with backward Euler as base method~\cite{QuEtAl2015}.
Figure~\ref{fig:stiff_limit} shows the spectral radius (left) and norm (right) of $\Tt{E}$ for the limit case~\eqref{eq:error_matrix_limit} and for~\eqref{eq:error_matrix} with a fast wave that is fifty or a hundred times faster than the slow wave.
Since the spectral radius remains smaller than unity up to large values of $M$ even for infinitely large $\lamfast$, \fwsw~still converges for $k \to \infty$.
For $M=12$, the spectral radius in the limit case finally becomes larger than unity (for $\lamfast = 100$, this happens for $M=11$) , but since $M=9$ e.g. would already allow to construct methods of order up to $17$ (using Radau nodes) this will most likely not be a relevant issue.
Note that, since the norm of $\Tt{E}$ is larger than one, convergence can be slow.
Modifications based on GMRES exist that can improve SDC convergence for stiff problems~\cite{HuangEtAl2006} but their exploration for \fwsw~is left for future work.

\subsection{Stability}\label{subsec:stability}
\begin{figure}[!th]
	\centering
	\includegraphics[scale=1]{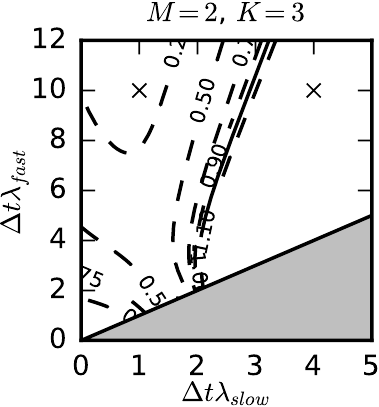}\hfill
	\includegraphics[scale=1]{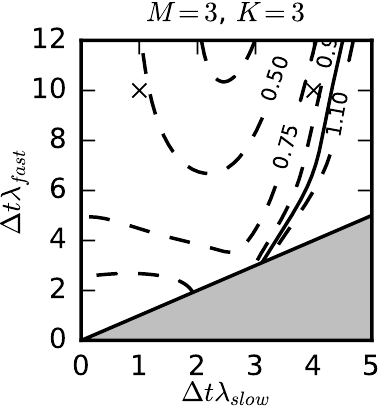}\hfill
	\includegraphics[scale=1]{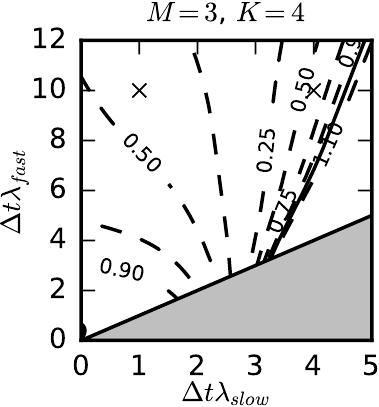}\\[1em]
	\includegraphics[scale=1]{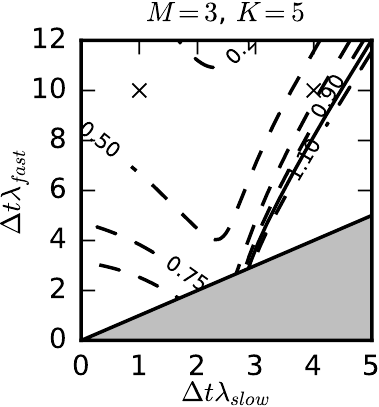}\hfill
	\includegraphics[scale=1]{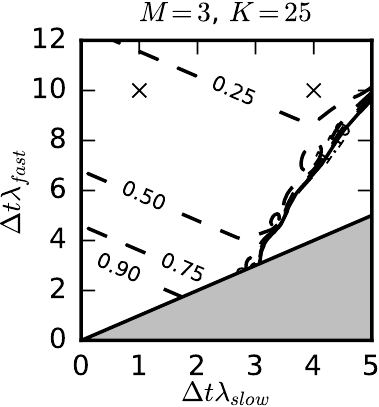}\hfill
	\includegraphics[scale=1]{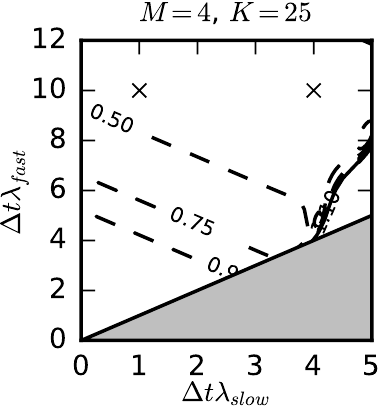}
	\caption{Stability domains of different configurations of \fwsw. $M$ indicates the number of quadrature nodes, $K$ the number of iterations. The gray region is where $\lamfast < \lamslow$ and the splitting becomes nonsensical. The values used for the plots in Figure~\ref{fig:stab_vs_k} are marked with crosses.}
	\label{fig:stability}
\end{figure}
Stability of SDC with splitting has been studied for the case where the fast dynamics correspond to negative real eigenvalues~\cite{LaytonMinion2005}.
First results on the stability of \fwsw~also exist~\cite{Weingartz2014}, but for Gauss-Lobatto nodes and without the derivation of a stability function.
Here, to study stability, we derive a formula for the update from $u_0$ to $u_{n+1}$.
Denote the left hand side matrix in~\eqref{eq:scalar-zero-to-node-system} as $\Tt{L}$  and the matrix on the right hand side as $\Tt{R}$ so that~\eqref{eq:scalar-zero-to-node-system} becomes 
\begin{equation}
\Tt{L} \Tt{U}^{k+1} = \Tt{U}_0 + \Tt{R} \Tt{U}^k,
\end{equation}
where both $\Tt{L}$ and $\Tt{R}$ depend on $\Delta t \lamfast$ and $\Delta t \lamslow$.
Using induction, it is straightforward to show that
\begin{equation}
	\Tt{U}^k = \left( \Tt{L}^{-1} \Tt{R} \right)^k \Tt{U}_0 + \sum_{j=0}^{k-1} \left( \Tt{L}^{-1} \Tt{R} \right)^j \Tt{L}^{-1} \Tt{U}_0.
\end{equation}
Denoting $\Tt{q} := \left( q_1, \ldots, q_M \right)$ and $\Tt{1} = (1, \ldots, 1)^{\text{t}}$, a full step of \fwsw~with $k$ iterations can be written as
\begin{equation}
	u_{n+1} = \left( 1 + i \lambda \Tt{q}  \left( \left( \Tt{L}^{-1} \Tt{R}\right)^k + \sum_{j=0}^{k-1} \left( \Tt{L}^{-1} \Tt{R} \right)^j \Tt{L}^{-1} \right) \Tt{1} \right) u_{0}
\end{equation}
so that the stability function of \fwsw~is given by
\begin{equation}
	\label{eq:stability_function}
	R(\Delta t \lamfast, \Delta t \lamslow) = 1 + i \lambda \Tt{q} \left( \left( \Tt{L}^{-1} \Tt{R}\right)^k + \sum_{j=0}^{k-1} \left( \Tt{L}^{-1} \Tt{R} \right)^j \Tt{L}^{-1} \right) \Tt{1}.
\end{equation}

Figure~\ref{fig:stability} shows the stability domains computed from~\eqref{eq:stability_function} for  different configurations of \fwsw~-- orders three, three and four in the upper row and orders five, five and seven in the lower.
Note that for the last two figures with $K=25$ the order is governed by the quadrature rule, not $K$.
The grey areas indicate $\lamfast < \lamslow$ where the splitting becomes nonsensical.

In all configurations, as long as $\Delta t \lamslow$ is small enough, the method remains stable for arbitrary large values of $\lamfast$.
While the y-axis in the figures goes only up to $\Delta t \lamfast = 12$, other experiments not documented here suggest that there is no stability limit on $\lamfast$. 
However, numerical damping becomes stronger as $\lamfast$ increases and $\left| R(\Delta t \lamfast, \Delta t \lamslow) \right|$ much smaller than unity.

In general, stability domains become larger when $K$ or $M$ is increased.
However, this does not happen monotonically and, in particular when increasing the number of iterations, the stability domain for $K+1$ does not always encompass the one for $K$.
For example, going from $K=4$ to $K=5$ for $M=3$ improves stability in some regions (upper right region) but slightly worsens it for smaller values of $\Delta t \lamfast$.
Similar behaviour is seen when increasing the number of quadrature nodes.
Going from $M=2$, $K=3$ to $M=3$, $K=3$ improves stability significantly, allowing for a slow CFL number of around three instead of two -- in this case, the stability domain for $M=3$ encompasses the one for $M=2$, but other examples can be found where this is not the case.
As $K \to \infty$, if SDC converges, it reproduces the stability properties of the underlying collocation method. 
The Radau based collocation method is stable everywhere so that instability regions of \fwsw~for $K=25$ indicate regions where the SDC iteration is not converging.
Note that the stability regions for $K=4$ and $K=5$ already match the eventual limit quite closely.

\begin{figure}[!t]
	\centering
	\includegraphics[scale=1]{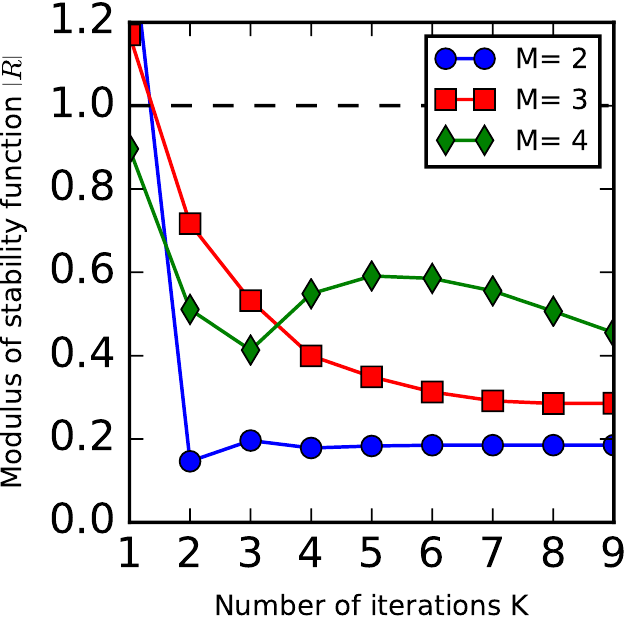}\hfill
	\includegraphics[scale=1]{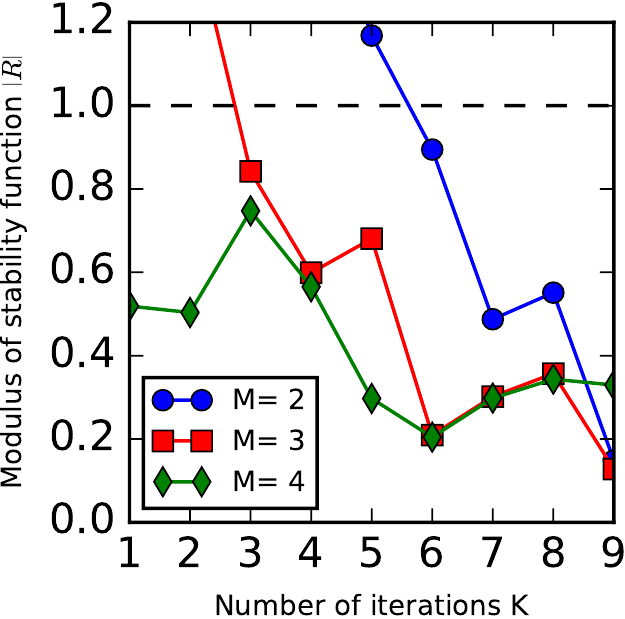}
	\caption{Modulus of the stability function for $\lamfast=10$ and $\lamslow=1$ (left) and $\lamslow=4$ (right) for different values of $M$.}
	\label{fig:stab_vs_k}
\end{figure}
Figure~\ref{fig:stab_vs_k} shows the modulus of the stability function versus $K$ for two fixed values of $\lamslow$, $\lamslow$; the points are marked with crosses in Figure~\ref{fig:stability}.
It illustrates again how larger values for $M$ and $K$ typically lead to better stability: in the left figure, for $\lamslow=1$, $M=4$ is stable for all values of $K$ while $M=2$ and $M=3$ are unstable for $K=1$ but stable for $K \geq 2$.
The influence of $M$ is more pronounced in the right figure where $\lamslow = 4$. For $M=2$, it takes six iterations for the method to become stable, for $M=3$ it still takes $K=3$ while for $M=4$ the method is stable throughout.
Note that the type of quadrature nodes can have a significant influence, see Section~\ref{subsec:quad_node_choice}.
\subsection{Dispersion relation}\label{subsec:dispersion}
To analyse the wave propagation characteristics of \fwsw~we derive the semi-discrete dispersion relation for the acoustic-advection equations
\begin{subequations}
\label{eq:acoustic_advection}
\begin{align}
	u_t + U u_x + c_s p_x &= 0  \\
	p_t + U p_x + c_s u_x &= 0
\end{align}
\end{subequations}
with a sound velocity $c_s$ that is significantly faster than the advection velocity $U$.
First, rewrite the system in matrix form
\begin{equation}
	\label{eq:acoustic_advection_system}
	\begin{pmatrix} u \\ p \end{pmatrix}_t = - \begin{pmatrix} U & 0 \\ 0 & U \end{pmatrix} \begin{pmatrix} u \\ p \end{pmatrix}_x  - \begin{pmatrix} 0 & c_s \\ c_s & 0 \end{pmatrix} \begin{pmatrix} u \\ p \end{pmatrix}_x.  
\end{equation}
The term with $U$ is treated explicitly, the acoustic term with $c_s$ implicitly.
Now assume a plane wave solution in space
\begin{equation}
	\label{eq:plane_wave_cont}
	u(x,t) = \hat{u}(t) e^{i \kappa x}, \quad p(x,t) =  \hat{p}(t) e^{i \kappa x}
\end{equation}
with wave number $\kappa$ so that~\eqref{eq:acoustic_advection_system} becomes
\begin{equation}
	\label{eq:acoustic_advection_system_fourier}
	\begin{pmatrix} \hat{u} \\ \hat{p} \end{pmatrix}_t =  -\Tt{U}_{\text{adv}} \begin{pmatrix} \hat{u} \\ \hat{p} \end{pmatrix} - \Tt{C}_s \begin{pmatrix} \hat{u} \\ \hat{p} \end{pmatrix}
\end{equation}
with
\begin{equation}
	\label{eq:CUmatrixDef}
	\Tt{C}_s := i \kappa \begin{pmatrix} 0 & c_s \\ c_s & 0 \end{pmatrix}, \quad \Tt{U}_{\text{adv}} := i \kappa \begin{pmatrix} U & 0 \\ 0 & U \end{pmatrix}.
\end{equation}
To obtain the dispersion relation of the fully continuous problem assume also a plane wave solution in time, that is
\begin{equation}
	\label{eq:plane_wave_time}
	\hat{u}(t) = u_0 e^{-i \omega t}, \quad \hat{p}(t) = p_0 e^{-i \omega t},
\end{equation}
with frequency $\omega$ so that~\eqref{eq:acoustic_advection_system_fourier} becomes
\begin{equation}
	\begin{pmatrix} -i \omega + i \kappa U & i \kappa c_s \\ i \kappa c_s & - i \omega + i \kappa U \end{pmatrix} \begin{pmatrix} u_0 \\ p_0 \end{pmatrix} = 0.
\end{equation}
For this system to have a solution for general values of $u_0$, $p_0$, the determinant of the matrix has to be zero, which gives the continuous dispersion relation of~\eqref{eq:acoustic_advection}
\begin{equation}
	\label{eq:dispersion_cont}
	\omega_{1,2} =  \left( U \pm c_s \right) \kappa.
\end{equation} 

To derive the semi-discrete dispersion relation of \fwsw~we apply it to~\eqref{eq:acoustic_advection_system_fourier}.
Since the problem now has two components, $u_0$ and $p_0$, the ``zero-to-node" SDC sweep~\eqref{eq:scalar-zero-to-node-system} for~\eqref{eq:acoustic_advection_system_fourier} becomes
\begin{align}
	&\left( \Tt{I} - \Delta t \left( \Qfast \otimes \Tt{C}_s \right) + \Delta t \left( \Qslow \otimes \Tt{U}_{\text{adv}} \right) \right) \Tt{X}^{k+1} \\ \nonumber
	&= \Tt{X}_0 - \Delta t \left( \Qfast \otimes \Tt{C}_s + \Qslow \otimes \Tt{U}_{\text{adv}} \right) \Tt{X}^k  
	+ \Delta t \Tt{Q} \otimes \left( \Tt{C}_s + \Tt{U}_{\text{adv}} \right) \Tt{X}^k
\end{align}
with
\begin{equation}
	\Tt{X} := \left( u_1, p_1, \ldots, u_M, p_M \right)^{\text{t}}, \quad \Tt{X}_0 := \left( u_0, p_0, u_0, p_0, \ldots, u_0, p_0 \right)^{\text{t}}.
\end{equation}
Here, the matrices $\Tt{U}_{\text{adv}}$ and $\Tt{C}_s$ essentially take the role of $\lamslow$ and $\lamfast$.
Therefore, equation~\eqref{eq:stability_function} for the stability function remains valid but with
\begin{equation}
	\Tt{L} := \left( \Tt{I} - \Delta t \left( \Qfast \otimes \Tt{C}_s + \Qslow \otimes \Tt{U}_{\text{adv}} \right) \right)
\end{equation}
and
\begin{equation}
	\Tt{R} := \Delta t \left( \Tt{Q} \otimes \left( \Tt{C}_s + \Tt{U}_{\text{adv}} \right) - \left(\Qfast \otimes \Tt{C}_s + \Qslow \otimes \Tt{U}_{\text{adv}}\right) \right),
\end{equation}
leading to the update formula
\begin{equation}
	\label{eq:sdc_acoustic_advection_update}
	\Tt{X}_{n+1} = \Tt{X}_0 + \left( \Tt{q} \otimes \left( \Tt{C}_s + \Tt{U}_{\text{adv}} \right) \right) \left( \left( \Tt{L}^{-1} \Tt{R}\right)^k + \sum_{j=0}^{k-1} \left( \Tt{L}^{-1} \Tt{R} \right)^j \Tt{L}^{-1} \right) \Tt{X}_0
\end{equation}
with $\Tt{X}_0 = \Tt{e} \otimes (u_0, p_0)$ and $\Tt{e} = (1,\ldots,1) \in \mathbb{R}^M$.
Now, instead of a continuous plane wave~\eqref{eq:plane_wave_time}, consider a solution in time of the form
\begin{equation}
	\label{eq:plane_wave_semi_disc}
	\hat{u}^{n} = u_0 e^{-i \omega n \Delta t}, \quad \hat{p}^n = p_0 e^{-i \omega n \Delta t},
\end{equation}
where $\hat{u}^n \approx \hat{u}(t_n)$, $\hat{p}^n \approx \hat{p}(t_n)$ are approximate solutions at some time step $t_n = n \Delta t$.
For a time stepping scheme with an update matrix $\Tt{Z}$, that is
\begin{equation}
	\begin{pmatrix} u \\ p \end{pmatrix}^{n+1} = \Tt{Z} \begin{pmatrix} u \\ p \end{pmatrix}^n,
\end{equation}
this ansatz gives
\begin{equation}
	\label{eq:semi_disc}
	\left[ \begin{pmatrix} e^{-i \omega \Delta t} & 0 \\ 0 & e^{-i \omega \Delta t} \end{pmatrix} - \Tt{Z} \right] \begin{pmatrix} u \\ p \end{pmatrix}^n = 0.
\end{equation}
Note that $\Tt{Z}$ does depend on $\Tt{U}_{\text{adv}}$ as well as $\Tt{C}_s$ and thus on $U$, $c_s$ and $\kappa$.
For \fwsw, the matrix $\Tt{Z}$ can be constructed by evaluating~\eqref{eq:sdc_acoustic_advection_update} for $(u_0,p_0)=(1,0)$ and $(u_0,p_0)=(0,1)$.
As in the continuous case, the dispersion relation corresponds to the roots of the determinant of the matrix in~\eqref{eq:semi_disc}. 
To compute the frequencies $\omega$ for a given wave number $\kappa$, the following equation has to be solved:
\begin{equation}
	\label{eq:dispersion_relation}
	\left( e^{-i\omega \Delta t} - \Tt{Z}_{11} \right) \left( e^{-i\omega \Delta t} - \Tt{Z}_{22} \right) - \Tt{Z}_{12}\Tt{Z}_{21} = 0
\end{equation}
where $\Tt{Z}_{11}$, $\Tt{Z}_{22}$, $\Tt{Z}_{21}$ and $\Tt{Z}_{12}$ are the entries of the matrix $\Tt{Z} \in \mathbb{C}^{2\times2}$.
We solve~\eqref{eq:dispersion_relation} using the symbolic Python package \emph{sympy}~\cite{sympy}.
\begin{remark}
To analyse dispersion when also the spatial derivative is discretised, assume a spatial solution of the form $e^{i \kappa \Delta x j}$ and replace the factor $i\kappa$ in~\eqref{eq:CUmatrixDef} with the symbol of a finite difference stencil, e.g. $\sin(\kappa \Delta x)/\Delta x$ for second-order centred differences~\cite[Sect. 3.3.1]{Durran2010}.
\end{remark}
\begin{figure}[!thp]
	\centering
	\includegraphics[scale=1]{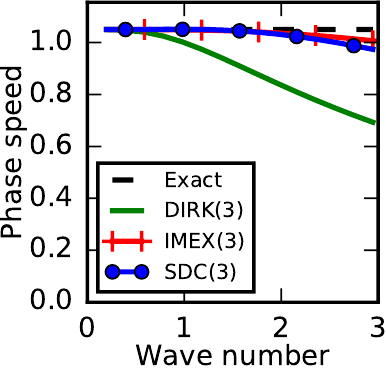}\hfill
	\includegraphics[scale=1]{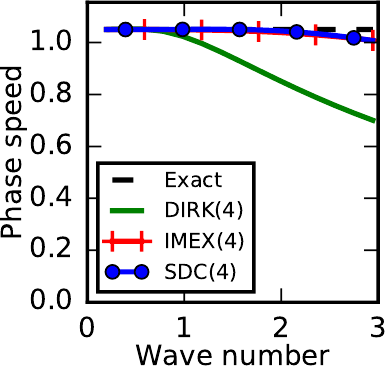}\hfill
	\includegraphics[scale=1]{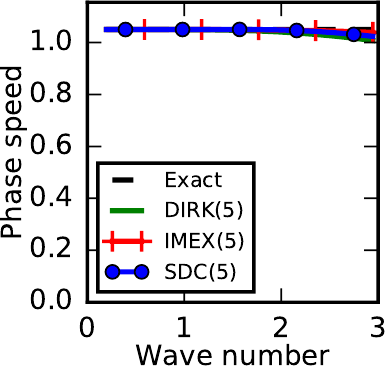}\\[1.0em]
	\includegraphics[scale=1]{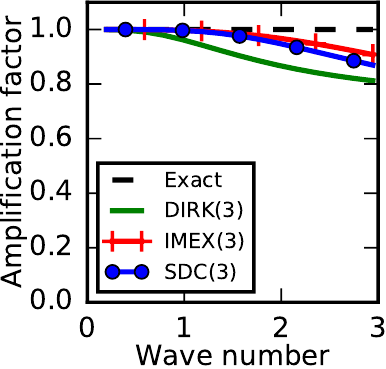}\hfill
	\includegraphics[scale=1]{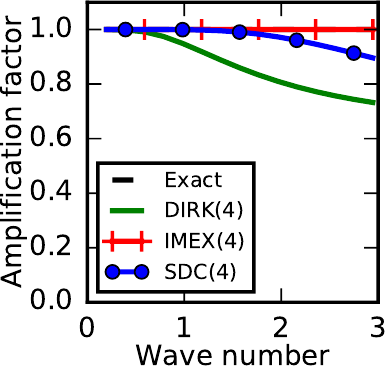}\hfill
	\includegraphics[scale=1]{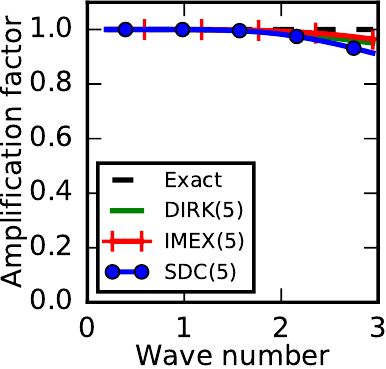}
	\caption{Semi-discrete dispersion relation for $U=0.05$ and $c_s=1.0$ for \fwsw, IMEX and DIRK methods of order three, four and five. Shown is the phase speed (upper) and amplification factor (lower) depending on the wave number $\kappa$.}		
	\label{fig:dispersion}
\end{figure}

Figure~\ref{fig:dispersion} shows the semi-discrete phase speed $\text{Real}(\omega)/\kappa$ and the amplification factor $\exp(\text{Imag}(\omega))$ for \fwsw, DIRK and IMEX methods of order three, four and five.

For order three, all three methods artificially slow down high wave number modes, but the effect is significantly more pronounced for DIRK(3) than for SDC(3) and IMEX(3).
All methods cause some attenuation particularly of high wave number modes, but again the effect is much more pronounced for DIRK(3) than for IMEX(3) and SDC(3).
The here presented variant of SDC uses $M=3$ nodes and $K=3$ iterations to achieve order three.
Interestingly, despite being formally of the same order of accuracy, third-order SDC with $M=2$ and $K=3$ (not shown) produces significantly stronger artificial slowing and damping.

For order four, phase speeds are almost identical to the exact values for IMEX and SDC except for minimal slowing of very large wave number modes.
In contrast, DIRK(4) does not provide a significant improvement compared to DIRK(3) and still produces inaccurate phase speeds across most of the spectrum.
In terms of dissipation, fourth-order \fwsw~produces slightly less artificial damping for very high wave number modes than SDC(3).
DIRK(4) shows significant attenuation across most of the wave number spectrum while IMEX(4) shows no numerical diffusion at all.

Lastly, all fifth-order methods give a quite accurate representation of the wave propagation characteristics of the continuous problem: there is very little slowdown and damping and only for high wave number modes.
Such semi-discrete propagation characteristics are attractive, because even for high frequency waves there are almost no phase speed errors and thus little numerical dispersion.
Also, low and medium wave number waves are propagated without amplitude errors while high wave number modes are slightly damped. 
While excessive numerical diffusion causes inaccurate solutions, a complete lack of numerical diffusion for large wave number modes retains spatially poorly resolved modes and can be problematic in atmospheric models with complex sub-scale models~\cite{VaterEtAl2011}.

\section{Numerical examples}\label{sec:numerics}
To demonstrate \fwsw's performance, numerical examples are presented below for a linear one-dimensional acoustic-advection problem with multi-scale initial data and for the two-dimensional compressible Boussinesq equations.

\subsection{Acoustic-advection}
To verify that \fwsw~provides the expected convergence order, consider the one-dimensional acoustic-advection problem~\eqref{eq:acoustic_advection} on a periodic domain $[0,1]$.
We split the equation according to
\begin{equation}
	\ffast(u,p) = \begin{pmatrix} c_s p_x \\ c_s u_x \end{pmatrix} \quad \text{and} \quad \fslow(u,p) = \begin{pmatrix} U u_x \\ U p_x \end{pmatrix}
\end{equation}
so that advection is treated explicitly while acoustic waves are integrated implicitly.
For initial data $u(x,0) \equiv 0 $ and $p(x,0) = p_0(x)$ the analytical solution of~\eqref{eq:acoustic_advection} reads
\begin{subequations}
\begin{align}
	u(x,t) &= \frac{1}{2} p_{0}\left( x - \left[U + c_{s} \right] t \right) - \frac{1}{2} p_{0}\left( x - \left[ U - c_{s} \right] t \right) \\
	p(x,t) &= \frac{1}{2} p_{0}\left( x - \left[U + c_{s} \right] t \right) + \frac{1}{2} p_{0}\left( x - \left[ U - c_{s} \right] t \right).
\end{align}
\end{subequations}
In line with the continuous dispersion relation~\eqref{eq:dispersion_cont} the solution consists of two modes travelling with phase velocities $c_{1,2} = \omega_{1,2}/\kappa = U \pm c_s$.
We set $T = 1.0$, $U = 0.1$ and $c_s = 1.0$.
The advective derivative is discretised with a fifth-order, the acoustic derivative with a sixth-order finite difference stencil.
All runs use five times as many spatial nodes as there are time steps, resulting in $C_{\text{fast}} = 5.0$ and $C_{\text{slow}} = 0.5$ in all runs, so that the fast mode is far from being well resolved.
Three configurations of \fwsw~are tested, all of them using $M=3$ Gauss-Radau nodes. 
The order is set by performing either $K=3$, $K=4$ or $K=5$ sweeps.

\begin{figure}[!th]
	\centering
	\includegraphics[scale=0.975]{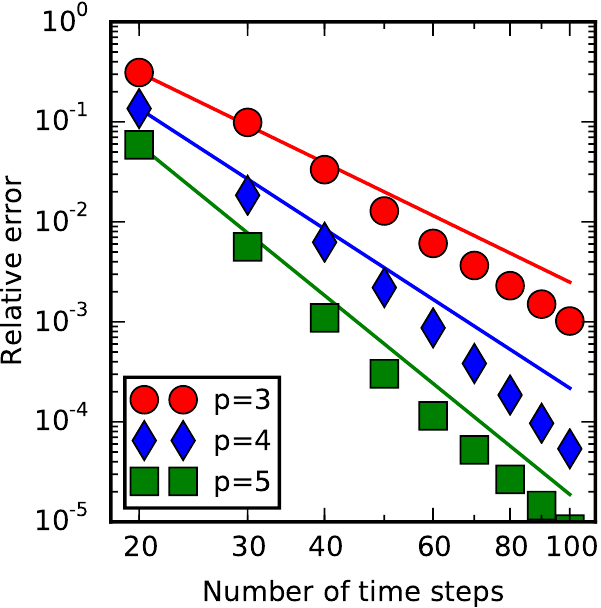}\hfill
	\includegraphics[scale=0.975]{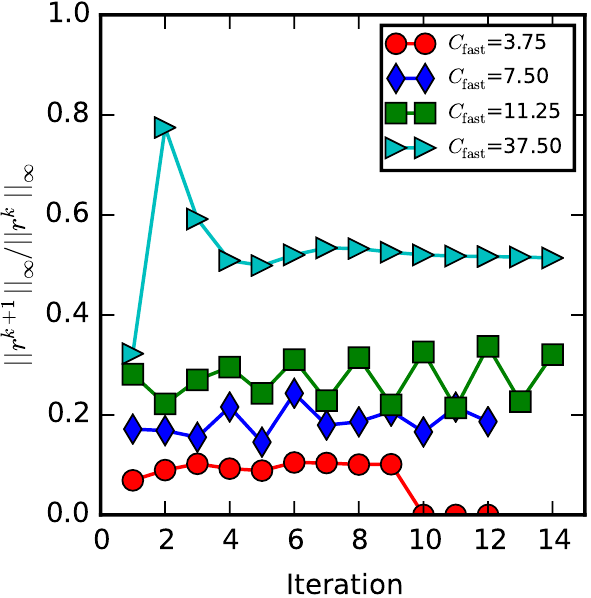}
	\caption{Left: Convergence of \fwsw~with orders three, four and five versus number of time steps. Both axes are scaled logarithmically. Right: Convergence rate of the \fwsw~iteration~for fixed $\Delta t$ and $\lamslow$ and varying values for $\lamfast$ versus the number of iterations $k$.}
	\label{fig:convergence}
\end{figure}

Figure~\ref{fig:convergence} (left) shows the relative error in the $\left\| \cdot \right\|_{\infty}$-norm at the end of the simulation, plotted against the number of time steps for $p_0(x) =  \sin(2\pi x) + \sin(5 \pi x)$.
As a guide to the eye, lines corresponding to orders three, four and five are drawn.
All three configurations of \fwsw~show the expected (or slightly better) order of convergence.
This illustrates that while the theoretical estimate of the convergence order shown above required $\Delta t\left| \lamfast \right| < 1$, in practice the expected order is observed much earlier.

In addition, the right graphic in Figure~\ref{fig:convergence} shows the ratio of SDC residuals from one sweep to the next for $M=3$ nodes over $15$ iterations.
The plotted ratio between residuals gives an estimate of the rate of convergence.
Here, a single time step of length $\Delta t = 0.025$ with $N_x = 300$ spatial nodes is performed for an advection velocity of $U=0.1$, corresponding to an advective CFL number of $C_{\text{slow}} = 0.75$.
Residuals are shown for four different values of sound speed $c_s$, leading to fast CFL numbers between $C_{\text{fast}} = 3.75$ and $C_{\text{fast}} = 37.5$.
For a large CFL number of $C_{\text{fast}} = 11.25$, \fwsw~still converges quickly with rates around $0.3$.
Even for an unrealistically large value of $C_{\text{fast}}=37.5$ \fwsw~still converges reasonably fast. 
Residuals are reduced in most iterations by a factor of about one half.
However, experiments not documented here suggest that if the fast wave speed is very large, much smaller time steps are needed to recover the expected order of convergence in $\Delta t$.

\subsection{Acoustic-advection with multi-scale initial data}\label{subsec:multiscale}
To assess how well \fwsw~damps highly oscillatory modes, we study an example from Vater et al.~\cite{VaterEtAl2011} with multi-scale initial data.
Let
\begin{equation}
	p(x,0) = p_0(x - x_0) + p_1(x - x_1)
\end{equation}
and $u(x,0) = p(x,0)$.
This results in a purely rightward travelling solution.
In contrast to Vater et al., we use a non-zero advection velocity $U=0.05$ and also a non-staggered mesh.
The purely large scale initial data is given
\begin{equation}
	p_0(x) = \exp\left( -\frac{x^2}{\sigma_0^2} \right)
\end{equation}
with $x_0 = 0.75$, $\sigma_0 = 0.1$ and $p_1 \equiv 0$.
The multi-scale initial data uses
\begin{equation}
	p_1(x) = p_0(x) \cos(k x / \sigma_0)
\end{equation}
with $x_1 = 0.25$ and $k = 7.2 \pi$ instead.
The domain is the unit interval $[0,1]$ with periodic boundary conditions and $N=512$ nodes in space.
The simulation is run until $T=3.0$ with $N_{\text{steps}} = 154$ time steps with $c_s = 1.0$, corresponding to an acoustic CFL number of $10$.
The advective CFL number is $0.5$.
\begin{figure}
	\centering
	\includegraphics[scale=1]{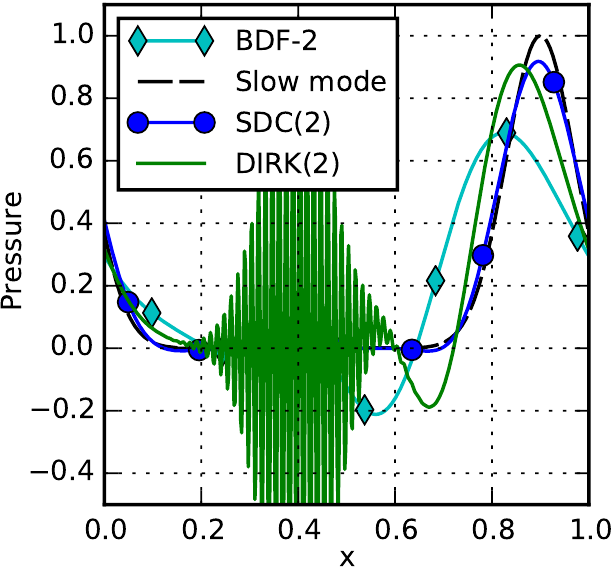}\hfill
	\includegraphics[scale=1]{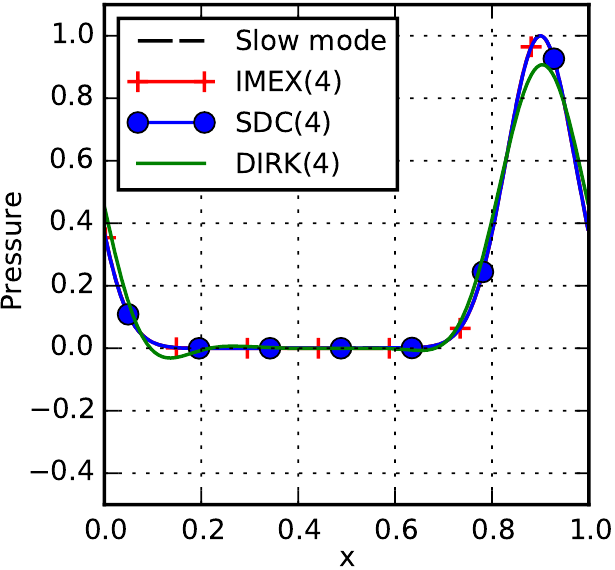}
	\caption{Numerical solution of the acoustic-advection equation with multi-scale initial data integrated with second order (left), using $M=2$, $K=2$ for SDC, and fourth order (right), using $M=3$, $K=4$ for SDC. Shown is the pressure $p$ at the final time $T=3$ when the slow part $p_0$ has been advected from $x_0 = 0.75$ to $x=0.9$ and the fast part $p_1$ has completed three revolutions. IMEX(2) is unstable and not plotted. The solutions provided by IMEX(4) and SDC(4) are indistinguishable in this plot.}
	\label{fig:sdc-multiscale}
\end{figure}

Figure~\ref{fig:sdc-multiscale} shows the solution produced by SDC, DIRK and IMEX methods of order two (left) and four (right). 
A backward differentiation formula (BDF) of order two is also run.
For comparison, the slow mode $p_0$ at the end of the simulation is plotted.
For SDC and DIRK, orders three and five (not shown) are similar to order four with somewhat more pronounced numerical diffusion for DIRK(3).
The IMEX methods of order two, three and five are unstable for this configuration.

Note that DIRK(2) corresponds to the midpoint rule which, for the linear problem studied here, is equivalent to the trapezoidal rule.
Both DIRK(2)/trapezoidal rule and BDF-2 match the results in Vater et al.:
BDF-2 removes the high frequency oscillations but introduces significant dispersion and also noticeable damping of the slow mode.
In contrast, DIRK(2) preserves the amplitude of the high frequency modes but slows them down to almost zero velocity.
Such undamped but wrongly propagated modes can have significant negative influence as discussed by Vater et al..
SDC(2) removes the high frequency waves, just as BDF-2, but also correctly propagates the slow mode without discernible dispersion and only little attenuation.

All three investigated fourth-order methods produce good solutions.
DIRK(4) shows some dispersion, in line with the too slow discrete phase speeds diagnosed in Section~\ref{subsec:dispersion}, and visible damping of the slow mode.
In contrast, both SDC(4) and IMEX(4) manage to damp the high frequency oscillations while still correctly advecting the slow mode without any discernible loss of amplitude.
Both solutions are indistinguishable in the plot.
\subsection{Compressible Boussinesq equations}
A key advantage of \fwsw~is that order of accuracy can be arbitrarily increased by simply adjusting run time parameters $K$ and $M$.
While the results so far suggest that \fwsw~provides more accurate solutions than its DIRK counterpart and solutions comparable to IMEX, it also requires significantly more evaluations of the right-hand side.
DIRK(4), for example, requires four (potentially nonlinear) implicit solves per time step, IMEX(4) requires six linear solves while fourth-order \fwsw~with $M=3$ and $K=4$ requires twelve.
However, for PDEs, the cost of each of these solves is not constant but depends on the number of iterations required by the employed solver.
The iterative nature of SDC provides increasingly accurate initial guesses which can reduce the cost of later sweeps~\cite{SpeckEtAl2015_DDM}.
We demonstrate that \fwsw~can outperform DIRK and compete with IMEX.

As the second and more complex test problem, we study the linearised Boussinesq equations governing compressible flow of a stably stratified fluid
\begin{subequations}
\begin{align}
	u_t + U u_x + p_x &= 0 \\
	w_t + U w_x + p_z &= b \\
	b_t + U b_x + N^2 w &= 0 \\
	p_t + U p_x + c_s^2 \left( u_x + w_z \right) &= 0.
\end{align}
\end{subequations}
They can be derived from the linearised Euler equations by a transformation of variables~\cite[Section 8.2]{Durran2010}.
This system supports gravity and acoustic waves as well as advective motion due to the background velocity $U$.
For SDC and IMEX we split the equations as
\begin{equation}
	\ffast(u,w,b,p) = \begin{pmatrix} -p_x \\ b - p_z \\ - N^2 w \\ -c_s^2 \left( u_x + w_z \right) \end{pmatrix} \quad \text{and} \quad \fslow(u,w,b,p) = -U \begin{pmatrix}  u_x \\  w_x \\  b_x  \\  p_x \end{pmatrix},
\end{equation}
so that terms corresponding to acoustic and gravity waves are integrated implicitly while the slow advection is treated explicitly.
The DIRK method treats both terms implicitly.

We choose a standard configuration where a non-hydrostatic gravity wave propagates through a channel of length $\SI{300}{\kilo\metre}$ and height $\SI{10}{\kilo\metre}$~\cite{SkamarockKlemp1994}.
Velocities $u$ and $w$ as well as pressure are set to zero initially.
An initial buoyancy perturbation
\begin{equation}
	b(x,z,0) = d\theta \frac{ \sin\left(  \frac{\pi z}{H} \right) }{1 + \left( x - x_0 \right)^2/a^2}
\end{equation}
with $d\theta = 0.01$, $H=\SI{10}{\kilo\metre}$, $x_0 = \SI{50}{\kilo\metre}$ and $a=\SI{5}{\kilo\metre}$ is placed at $x = -\SI{50}{\kilo\metre}$, which generates waves propagating to both sides.
Periodic boundary conditions in the horizontal and no-slip boundary conditions at the top and bottom are employed.
Fifth-order upwind finite differences are used to discretise the advective derivatives and fourth-order centred differences for the acoustic derivatives.
\begin{figure}[!t]
	\centering
	\includegraphics[scale=1]{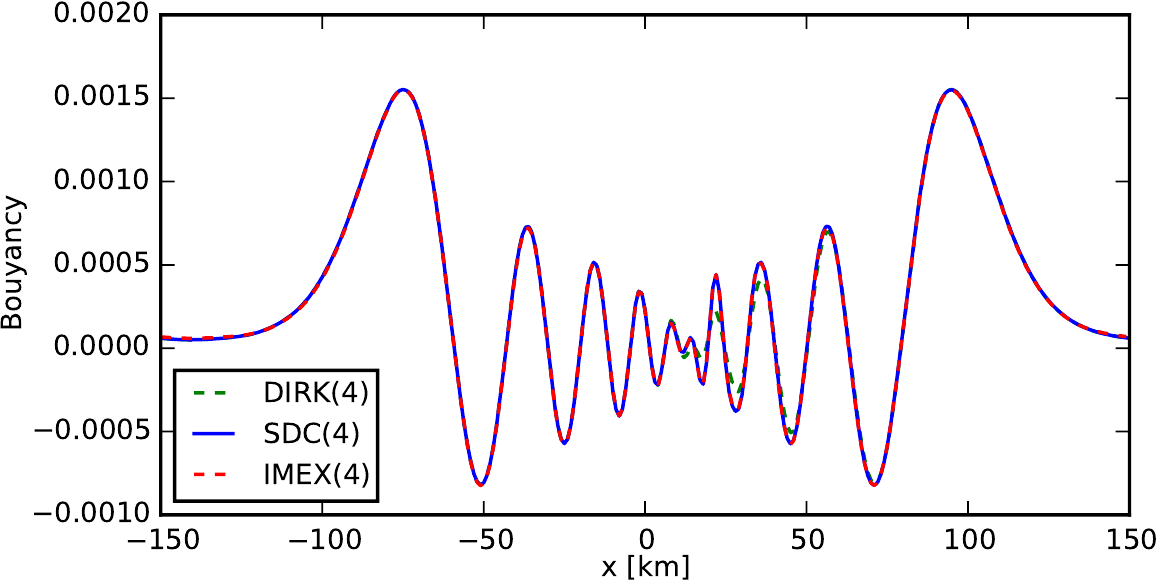}
	\caption{Cross section of the buoyancy $b$ at $z = \SI{5}{\kilo\metre}$ at $T=\SI{3000}{\second}$, computed with fourth-order \fwsw~as, DIRK and IMEX and $\Delta t = \SI{30}{\second}$. The solution from SDC(4) and IMEX(4) are indistinguishable.}
	\label{fig:boussinesq}
\end{figure}

The spatial resolution is $300 \times 30$ nodes, corresponding to $\Delta x = \SI{1}{\kilo\metre}$ and $\Delta z = \SI{0.32}{\kilo\metre}$.
The advection velocity is set to $U = \SI{20}{\metre\per\second}$, the acoustic velocity to $c_s = \SI{300}{\metre\per\second}$ and the stability frequency to $N = \SI{0.01}{\per\second}$.
We run the simulation until $T=\SI{3000}{\second}$ with a time step of either $\Delta t = \SI{30}{\second}$ or $\Delta t = \SI{6}{\second}$.
For the large time step, the resulting advective CFL number is $0.6$, the horizontal acoustic CFL number is $9.0$ while the vertical acoustic CFL number is $27.9$.
For the small time step, they are $0.12$, $1.80$ and $5.58$.
To solve the linear systems arising in the DIRK method and the implicit parts of \fwsw~and IMEX, the GMRES solver of the \emph{SciPy} package~\cite{scipy2001} is used with a tolerance of $10^{-5}$ and restart after $10$ iterations (the default values).
For SDC, to avoid over-solving in early sweeps, a tolerance equal to a factor times the SDC residual or the default is used, whatever is higher.
The factor is set to $0.1$ for all runs.
To estimate the temporal discretisation error, a reference solution is computed using fifth-order IMEX with a ten times smaller time step and a GMRES tolerance of $10^{-10}$.
Variants of each method of orders three, four and five are run and the final error is estimated against the reference solution.
Also, the total number of required GMRES iterations is logged.
SDC uses $M=3$ nodes with $K=3$, $K=4$ and $K=5$  iterations to realise the different orders.
\begin{table}[!t]
\centering
\begin{tabular}{|c|ccc|ccc|} \hline
\textbf{Third-order}     & \multicolumn{3}{c|}{$\Delta t = \SI{30}{\second}$} & \multicolumn{3}{c|}{$\Delta t = \SI{6}{\second}$}  \\
                                   & DIRK           & IMEX     & SDC           & DIRK            & IMEX         & SDC \\[0.2em] \hline
\# implicit solves         & 200              &               &  900           &  1000           &  2000         &  4500   \\
\# GMRES iterations  & \SI{46702}{} &               & \SI{25819} & \SI{28863}{} & \SI{13782} & \SI{25051}{}   \\
avg. it. per call            & 233.5           &               & 28.7          & 28.9              & 6.9             & 5.6          \\  
est. error                     & 1.8e-1         &  unstable &    1.1e-1   & 9.6e-2          &  1.7e-2       & 1.5e-2  \\ \hline       
\end{tabular}\vspace*{1ex}
\begin{tabular}{|c|ccc|ccc|} \hline
\textbf{Fourth-order}     & \multicolumn{3}{c|}{$\Delta t = \SI{30}{\second}$} & \multicolumn{3}{c|}{$\Delta t = \SI{6}{\second}$}  \\
                                   & DIRK              & IMEX        & SDC           & DIRK            & IMEX         & SDC \\[0.2em] \hline
\# implicit solves         & 300                & 500            &  1200         & 1500             & 2500         & 6000           \\
\# GMRES iterations  & \SI{100651}{} &  \SI{38092} & \SI{31105} & \SI{66136}{} & \SI{24068} & \SI{32696}{}   \\
avg. it. per call            & 335.5             &  76.2          & 25.9           & 44.1              & 9.6            & 5.4            \\  
est. error                     & 1.5e-1            &    1.3e-1    & 9.9e-2         & 9.4e-2          & 4.2e-3       & 2.9e-3         \\ \hline     
\end{tabular}\vspace*{1ex}
\begin{tabular}{|c|ccc|ccc|} \hline
\textbf{Fifth-order}     & \multicolumn{3}{c|}{$\Delta t = \SI{30}{\second}$} & \multicolumn{3}{c|}{$\Delta t = \SI{6}{\second}$}  \\
                                   & DIRK           & IMEX      & SDC            & DIRK           & IMEX          & SDC \\[0.2em] \hline
\# implicit solves         & 500              &                & 1500           & 2500           & 3500           & 7500           \\
\# GMRES iterations  & \SI{38334}{} &                &  \SI{34732} & \SI{24592}{} & \SI{24649} & \SI{32724}{}   \\
avg. it. per call            & 76.7             &                &  23.2           & 9.8               & 7.0             & 4.4            \\  
est. error                     & 9.6e-2          & unstable &  9.7e-2        & 3.4e-3          & 2.7e-3       &  2.6e-3        \\ \hline        
\end{tabular}\vspace*{1ex}
\caption{Number of implicit solves and total number of required GMRES iterations for the solution of the Boussinesq equations for DIRK, IMEX and \fwsw~of orders three, four and five.}
\label{table:gmres}
\end{table}

Figure~\ref{fig:boussinesq} shows a cross section through the buoyancy field $b$ at a height $z = \SI{5}{\kilo\metre}$ at the end of the simulation.
Gravity waves are propagating to the left and right and advection has moved the centre point by $\SI{60}{\kilo\metre}$ to the right, from $x = -\SI{50}{\kilo\metre}$ to $x=\SI{10}{\kilo\metre}$.
All methods properly resolve the larger scale oscillations at the fronts of the wave train.
For the small scale oscillations in the centre, DIRK(4) produces wave positions in line with SDC and IMEX but with slightly damped amplitudes.

Table~\ref{table:gmres} shows the total number of implicit solves over the course of the simulation, total number of required GMRES iterations, the average number of iterations per solve and the estimated error.
For order three, SDC(3) and DIRK(3) are stable for the large time step while IMEX is unstable.
SDC(3) is more accurate than DIRK(3) and requires significantly fewer GMRES iterations.
Interestingly, the third-order version of SDC using only $M=2$ nodes (not shown) requires \emph{more} overall GMRES iterations than for $M=3$ (\SI{29337} versus \SI{25819}), even though it requires only six solves per time step for a total of 600.
For the small time step, all methods are stable.
DIRK(3) is the most expensive, IMEX(3) the cheapest and SDC(3) in the middle.
SDC(3) is the most accurate method, but IMEX is comparable.

For the fourth-order methods with large time step, SDC is the cheapest and most accurate of the three methods.
When the time step is decreased, IMEX becomes the cheapest method, but SDC remains the most accurate.
In all configurations, SDC requires the fewest iterations per solve.
Note that when spatial resolution is increased and the system to be solved becomes larger, the number of GMRES iterations increases for all methods but the ordering seems to be unaffected.

Finally, for fifth-order with large time step, IMEX is unstable while both DIRK and SDC generate roughly the same error with SDC being about 10\% cheaper.
For the smaller time step, DIRK and IMEX are comparable in the number of required GMRES iterations with IMEX being more accurate.
SDC is more costly but slightly more accurate than IMEX.

These results are preliminary and a detailed, fair comparison of all three methods would probably warrant a paper on its own.
In particular, only a single problem and neither the effect of preconditioning the linear systems nor the influence of a nonlinear Newton solver are investigated here.
Nevertheless, these results illustrate that, despite the fact that it needs more implicit solves, SDC can be competitive compared to both DIRK and IMEX methods.
A more comprehensive comparison is planned for future work.

\subsection{A comment on the choice of quadrature nodes}\label{subsec:quad_node_choice}
For semi-implicit SDC applied to problems of advective-diffusive type, choosing Gauss-Lobatto nodes leads to good stability properties~\cite{LaytonMinion2005}.
We found this to be different for the fast-wave slow-wave case: when using the ``correct'' collocation update~\eqref{eq:picard_update}, stability regions are significantly smaller than for Radau or Legendre nodes (see also Remark~\ref{remark:update}).
In particular, Lobatto nodes lead to limits on $\Delta t \lamfast$ even for small values of $\Delta t \lamslow$.
Both Radau and Legendre nodes, in contrast, show good stability without a clear ranking: depending on the values for $M$ and $K$, one or the other can produce larger stability domains.
In terms of dispersion properties, Legendre and Radau nodes are comparable with Radau nodes causing slightly more numerical diffusion.
For the Boussinesq example, \fwsw~based on Radau nodes requires fewer overall GMRES iterations compared to Legendre nodes but the latter give slightly smaller errors.
In summary, all examples presented here were done using Gauss-Radau nodes but both types have advantages.
For the sake of brevity we do not present results for Legendre nodes but the interested reader could easily generate them using the published code~\cite{pySDC2016}.

\section{Conclusions}
The paper analyses semi-implicit spectral deferred corrections (SISDC) with fast-wave slow-wave splitting (\fwsw) where the stiff fast process is due to fast propagating waves instead of diffusion.
\fwsw~allows to easily construct splitting methods of arbitrary high order of accuracy.
The iteration error and local truncation error are analysed.
For the non-stiff limit, \fwsw~increases the order by one per iteration.
In the stiff limit, the error propagation matrix reduces to the non-split case with implicit Euler as base method.
Since the spectral radius remains smaller than unity, \fwsw~continues to converge but as the norm becomes larger than unity, convergence can become slow.
However, numerical examples suggest that even for rather large fast-wave CFL numbers, convergence is still reasonably good.
Stability function and semi-discrete dispersion relation are derived and analysed.
\fwsw~has good stability properties and phase and amplitude errors in line with Runge-Kutta IMEX methods of the same order.
Finally, performance is studied in numerical examples, showing that \fwsw~can be competitive with DIRK and IMEX methods in terms of cost and accuracy.

\section*{Acknowledgments}
All figures in this manuscript have been generated with the Python library \textit{matplotlib}~\cite{matplotlib}.
The source code used to generate the results in this paper is based on the Python framework \texttt{pySDC} and can be accessed through \emph{GitHub}~\cite{pySDC2016}.

\bibliographystyle{siam}
\bibliography{sdc,pint,refs,faulttol}

\end{document}